\documentclass[twoside]{amsart}
\usepackage{amsmath,amsthm,amssymb}

\theoremstyle{plain}
\newtheorem{theorem}{Theorem}[section]
\newtheorem{proposition}[theorem]{Proposition}
\newtheorem{lemma}[theorem]{Lemma}
\newtheorem{keylemma}[theorem]{Key-Lemma}

\theoremstyle{remark}
\newtheorem{remark}[theorem]{Remark}
\newtheorem{acknowledgments}{Acknowledgments}

\renewcommand{\k}{{\bf k}}
\newcommand{\tf}{\tilde f}
\newcommand{\dk}{\dot{k}}
\newcommand{\K}{{\bf K}}

\newcommand{\F}{{\mathbb F}}

\newcommand{\Z}{{\mathbb Z}}

\newcommand{\llceil}{\left\lceil}
\newcommand{\rrceil}{\right\rceil}
\newcommand{\llfloor}{\left\lfloor}
\newcommand{\rrfloor}{\right\rfloor}
\def\bar{\overline}

\newcommand{\HS}{{\mathcal{HS}}}

\newcommand{\ord}{ {\rm ord} }
\newcommand{\NP }{ {\rm  NP} }     
\begin{document}

\title{Hyperelliptic curves in characteristic 2}
\author{Jasper Scholten and Hui June Zhu}
\address{
Mathematisch Instituut,
Katholieke Universiteit Nijmegen,
Postbus 9010, 6500 GL Nijmegen.
The Netherlands.
}
\email{scholten@sci.kun.nl}

\address{
Department of mathematics,
University of California,
Berkeley, CA 94720-3840.
The United States.
}
\email{zhu@alum.calberkeley.org}

\date{\today}

\begin{abstract}
In this paper we prove that there are no hyperelliptic supersingular
curves over $\bar\F_2$ of genus $2^n-1$ for any integer $n\geq 2$.

Let $g$ be a natural number. Write $h=\llfloor\log_2(g+1)+1\rrfloor$.
Let $X$ be a hyperelliptic curve over $\bar\F_2$ of genus $g\geq 3$ of
$2$-rank zero, given by an affine equation $y^2-y=
\sum_{i=1}^{2g+1}c_ix^i$. We prove that the first
slope of the Newton polygon of $X$ is $\geq 1/h$. We also prove that
the equality holds if (I) $g<2^h-2$, $c_{2^h-1}\neq 0$; or (II)
$g=2^h-2$, $c_{2^h-1}\neq 0$ or $c_{3\cdot 2^{h-1}-1}\neq 0$.

Let $\HS_g/\bar\F_p$ be the intersection of
the supersingular locus with the open hyperelliptic Torelli locus in
the moduli space of principally polarized abelian varieties over
$\bar\F_p$ of dimensions $g$.  We show that $\dim\HS_g/\bar\F_p\leq
g-2$ for every $g\geq 3$.

We prove that $\dim\HS_4/\bar\F_2=2$ by showing that every
genus-4 hyperelliptic supersingular curve over $\bar\F_2$ has an
equation $y^2-y=x^9+c_5x^5+c_3x^3$ for some $c_5,c_3\in\bar\F_2$.
\end{abstract}

\maketitle

\section{Introduction}
        \label{section1}

Let $p$ be a prime number.  For any natural number $\nu$ let
$\F_{p^\nu}$ be a finite field of $p^\nu$ elements.  Let
$\bar\F_{p^\nu}$ be an algebraic closure of $\F_{p^\nu}$.  In this
paper a curve is a smooth, projective and geometrically integral
algebraic variety of dimension $1$.  A curve is {\em supersingular} if
its Jacobian is supersingular as an abelian variety (see
also~\cite{Oort:89}~\cite{Li-Oort:98}).
Newton polygon of a curve is the Newton polygon of
the Jacobian of the curve~\cite{Oort:91}.
Any curve of genus $1$ or $2$ and of $p$-rank zero has its
Newton polygon equal to a straight line segment of slope $1/2$. In
other words, it has to be supersingular.
The smallest possible first slope of the Newton polygon
of a  curve of genus $g\geq 2$ with $p$-rank zero is $1/g$.

In this paper, let $X$ be a hyperelliptic curve of genus $g$ over
$\bar\F_2$ given by an affine equation $y^2-y = {\tilde f}(x)$ where
$\tf(x)=\sum_{i=1}^{2g+1}c_{i}x^i$.  Denote by
$\NP_1(X)$ the first slope of the Newton polygon of of $X$.

\begin{theorem}\label{MT:3}
Let $g\geq 3$. Write $h=\llfloor\log_2(g+1)+1\rrfloor$.
Let $X$ be a hyperelliptic curve
over $\bar\F_2$ given by an equation $y^2 - y = \sum_{i=1}^{2g+1}c_ix^i$
where $c_{2g+1}=1$.  Then
\begin{enumerate}
\item[(I)] $\NP_1(X)\geq\frac{1}{h}$.
\item[(II)]
If $g<2^h-2$ and $c_{2^{h}-1}\neq 0$ then
$\NP_1(X)=\frac{1}{h}$.
\item[(III)]
If $g=2^h-2$ and
$c_{2^{h}-1}\neq 0$ or $c_{3\cdot2^{h-1}-1}\neq 0$ then
$\NP_1(X) = \frac{1}{h}$.
\end{enumerate}
\end{theorem}

Note that $g\leq 2^h-2$.  This theorem will be proved in
Section~\ref{section3}.  Note that every hyperelliptic curve over
$\bar\F_2$ of genus $g$ and of $2$-rank $0$ has an equation as in the
above theorem (see Proposition~\ref{P:1}). If $g=2^n-1$ for some
$n\geq 2$ then we have $g=2^{h-1}-1<2^h-2$ and
$c_{2^h-1}=c_{2g+1}=1\neq 0$.  Thus Theorem~\ref{MT:3} implies the
following theorem. It stands in striking comparison to the existence
of supersingular curves over $\bar\F_2$ of every genus~\cite{Geer:95}.

\begin{theorem}\label{MT:1}
For every integer $n\geq 2$, there are no hyperelliptic supersingular
curves over $\bar\F_2$ of genus $2^n-1$.
\end{theorem}

Let $\HS_g/\bar\F_p$ denote the intersection of the supersingular
locus with the open hyperelliptic Torelli locus in the moduli space of
principally polarized abelian varieties of dimension $g$ over
$\bar\F_p$ (see~\cite[Appendix, A.11]{Li-Oort:98}).
It is well-known that $\dim\HS_1/\bar\F_p=0$: For $p=2$,
there is exactly one
supersingular elliptic curve given by $y^2-y =x^3$.
For $p>2$, an elliptic
curve in Legendre form $y^2=x(x-1)(x-\lambda)$ is supersingular if and
only if $\lambda$ is a root of the Hasse polynomial as in Theorem
4.1~(b) in~\cite[V.4]{Silverman:86}.  There are finitely many of them.
One knows that $\HS_2/\bar\F_p$ is an open subset of the supersingular
locus, which is of dimension 1 (see~\cite{Li-Oort:98}).
In particular, a hyperelliptic
supersingular curve of genus $2$ over $\bar\F_2$
is of the form $y^2-y = x^5+c_3x^3$ for some $c_3\in\bar\F_2$
(see~\cite[Section~2]{Igusa:60}). It was proved that
$\HS_3/\bar\F_2$ is empty (see~\cite{Oort:89}). Higher genera cases
were considered in~\cite{Geer:92}~\cite{Geer:95}, where it was shown
that $\dim\HS_{2^n}/\bar\F_2 \geq n$ for every $n\geq 1$
and hyperelliptic curves over $\bar\F_2$
of the following form are supersingular:
\begin{eqnarray}\label{E:geer}
y^2-y &=& \sum_{i=0}^{n+1}c_{2^i+1}x^{2^i+1}.
\end{eqnarray}
See~\cite[Section~7]{SZ:2} for another proof.

There are also results on the intersection of the supersingular locus
and the closed hyperelliptic Torelli locus (see~\cite[Appendix,
A.11]{Li-Oort:98}).  Oort proved that for $g=3$ this intersection is
pure of dimension $1$ for all $p$ (see~\cite{Oort:89}).

\begin{theorem}\label{MT:4}
For $g\geq 3$ a generic hyperelliptic curve over $\bar\F_2$ of genus $g$ and
$2$-rank zero is not supersingular. Moreover,
$\dim\HS_g/\bar\F_2\leq g-2$. If $g=2^h-2$ for some $h\geq 3$ then
$\dim\HS_g/\bar\F_2\leq g-3$.
\end{theorem}

This upper bound of $\dim\HS_g/\bar\F_2$ can be achieved for $g=4$.

\begin{theorem}\label{MT:2}
If $g=4$ then $X$ is supersingular if and only if $X$ has an equation
$y^2-y = x^9+c_5x^5+c_3x^3$ for some $c_5,c_3\in\bar\F_2$.  In
particular, $\HS_4/\bar\F_2$ is irreducible of dimension $2$.
\end{theorem}

Theorems~\ref{MT:4} and~\ref{MT:2} will be proved in
Section~\ref{section4}.

\begin{acknowledgments}
The second author thanks Carel Faber for a conversation which
stimulated our computation and the research. It is a pleasure to
thank Shuhong Gao and Daqing Wan for helpful discussions. In
particular, we thank Bjorn Poonen for comments on an early version. A
part of this work was carried out while the second author was
associated with and supported by MSRI (at Berkeley).
\end{acknowledgments}

\section{$2$-adic box analysis}
           \label{section2}

Let $g\geq 3$.
For any nonnegative integer $m$ let
$(m)_2$ denote the binary expansion of $m$, and let $s(m)$
denote the sum of digits in $(m)_2$. Let $d=2g+1$ and
$h=\llfloor\log_2(g+1)+1\rrfloor$.  Every hyperelliptic curve $X$ over
$\bar\F_2$ of genus $g$ and 2-rank zero has an affine equation
\begin{eqnarray}\label{E:original}
y^2 - y = \sum_{i=1}^{d}c_ix^i
\end{eqnarray}
where $c_1,\cdots,c_{d}\in\F_q$ for some $2$-power $q$ and $c_{d}=1$.
(See Proposition~\ref{P:1}.)
Let $W(\F_q)$ be the ring of Witt vectors over $\F_q$,
and $f(x)=\sum_{i=1}^da_ix^i\in W(\F_q)[x]$ such that $a_i\equiv c_i\mod 2$
for $1\leq i\leq d$ and $a_d=1$. For nonnegative integers $N$ and $r$
define $C_r(N)$ by the following power series expansion:
\begin{eqnarray}\label{E:define-c}
(1+4f(x))^{(2^N-1)/2}
&=&\sum_{k_1=0}^{\infty}4^{k_1}\binom{(2^N-1)/2}{k_1}f(x)^{k_1}\\\nonumber
&=&\sum_{r=0}^{\infty}C_r(N)x^r.
\end{eqnarray}
For the ease of formulation, we define $C_r(\cdot)=0$ for all $r<0$.

The following lemma is borrowed from the Key-Lemma in~\cite{SZ:2}.
In that paper the curve $X/\F_q$ is lifted to a curve over $W(\F_q)$
given by the equation $y^2-y = f(x)$, and the coefficients of the
power series expansion of $(2y-1)^{2^N-1}$ in $x$ are considered.
Note that these are equal to the $C_r(N)$, since
$(2y-1)^2=1+4f(x)$.

\begin{keylemma}\label{keylemma}
Let $\lambda$ be a rational number with $0\leq \lambda\leq\frac{1}{2}$.
\begin{enumerate}
\item[i)] If for all $m\geq1$, $n\geq1$ and $1\leq j\leq g$,
\begin{eqnarray*}
\ord_2(C_{m2^{n+g-1}-j}(n+g-2))&\geq&\llceil n\lambda\rrceil
\end{eqnarray*}
then $$\NP_1(X)\geq\lambda.$$
\item[ii)] Let $1\leq j\leq g$ and $n_0\geq1$.
Suppose that for all
$m\geq1$ and $1\leq n<n_0$,
\begin{eqnarray*}
\ord_2(C_{m2^{n+g-1}-j}(n+g-2))&\geq&\llceil n\lambda\rrceil;
\end{eqnarray*}
suppose that for all $m\geq2$,
\begin{eqnarray*}
\ord_2(C_{m2^{n_0+g-1}-j}(n_0+g-2))&\geq&\llceil n_0\lambda\rrceil;
\end{eqnarray*}
suppose
\begin{eqnarray*}
\ord_2(C_{2^{n_0+g-1}-j}(n_0+g-2))&<& \llceil n_0\lambda\rrceil.
\end{eqnarray*}
Then
\begin{eqnarray*}
\NP_1(X)&<&\lambda.
\end{eqnarray*}
\end{enumerate}
\end{keylemma}

For any nonnegative integer $r$, let $\K_r$ be the subset of $\Z^d$
defined by
$$
\K_r:=\{{ }^t\!(k_1,k_2,\ldots,k_d)\in\Z^d\mid
k_1\geq k_2\geq\cdots\geq k_d\geq 0, \sum_{\ell=1}^{d} k_\ell = r
\}.
$$
For ease of formulation  we set $0^0:=1$, then we get by definition
\begin{eqnarray}
\label{E:bin}
C_r(N)=
\sum_{\k\in \K_r}4^{k_1}\binom{(2^N-1)/2}{k_1}
\prod_{\ell=1}^{d-1}\binom{k_{\ell}}{k_{\ell+1}}
\prod_{\ell=1}^{d-1}a_\ell^{k_{\ell}-k_{\ell+1}}.
\end{eqnarray} From the fact that $\ord_2(m!)=m-s(m)$ one can derive that
\begin{eqnarray*}
s(\k)
&:=& \ord_2\left(4^{k_1}\binom{(2^N-1)/2}{k_1}\prod_{\ell=1}^{d-1}
\binom{k_{\ell}}{k_{\ell+1}}\right)\\
&=& s(k_1)+\sum_{\ell=1}^{d-1}\ord_2(\binom{k_\ell}{k_{\ell+1}})\\
&=& s(k_1-k_2)+s(k_2-k_3)+\ldots + s(k_{d-1}-k_{d}) + s(k_{d}).
\end{eqnarray*}
Note that it is independent of $N$.
In fact, all properties of the coefficients $C_r(N)$ that we
need in this paper do not depend on $N$. For this reason, in the
rest of this paper we will drop the $N$ from the notation, and
simply write $C_r$ for these coefficients.

Let $\k=\ ^t\!(k_1,\ldots,k_d)\in\K_r$.
Let $k_d=\sum_{v\geq 0}k_{d,v}2^v$ be the binary expansion of
$k_\ell$, we introduce a {\em dot representation}
$$\dk_\ell:= [\ldots, \dk_{\ell,2},\dk_{\ell,1}, \dk_{\ell,0}]$$
in the following way: for $\ell=d$, let $\dk_{d,v}=k_{d,v}$
for all $v\geq 0$; for $1\leq \ell<d$, it is defined inductively by
$$\dk_{\ell-1,v}:=\dk_{\ell,v}+ \mbox{$2^v$-coefficient in the binary
expansion of $(k_{\ell-1}-k_{\ell})$},$$
for all $v\geq 0$.
It can be verified that
$k_\ell=\sum_{v\geq 0}\dk_{\ell,v}2^v$ for $1\leq \ell\leq d$.  Since
$k_\ell\geq k_{\ell+1}$ we have $\dk_{\ell-1,v}\geq \dk_{\ell,v}$ for
all $v$. It is not hard to observe
$$s(\k)=\sum_{v\geq 0}\left(
\sum_{\ell=1}^{d-1}(\dk_{\ell,v}-\dk_{\ell+1,v})+\dk_{d,v}\right)
= \sum_{v\geq 0}\dk_{1,v}.$$
We call the representation
${}^t\!(\dk_1, \dk_2,\ldots,\dk_d)$
{\em the 2-adic box} of $\k$, denoted by $\fbox{\k}$ for short. See below.
$$
\fbox{\k}:=
\begin{bmatrix}
\ldots &\dk_{1,2}&\dk_{1,1}&\dk_{1,0}\\
\ldots &\dk_{2,2}&\dk_{2,1}&\dk_{2,0}\\
&&\vdots&\\
\ldots &\dk_{d,2}&\dk_{d,1}&\dk_{d,0}
\end{bmatrix}.
$$

\begin{lemma}\label{miracle}
For $\k\in\K_r$ with $r\geq 1$, we have $s(\k)\geq
\llceil\frac{s(r)}{h}\rrceil.$ If $s(\k)=\frac{s(r)}{h}$ then
\begin{enumerate}
\item[1)] the $2$-adic box $\fbox\k$ consists of only $0$ and $1$'s;
\item[2)] for every $v\geq 0$
there are $\frac{s(r)}{h}$ many $\dk_{1,v}$'s equal to $1$ while
the others are $0$;
\item[3)] for every $v\geq 0$,
$s(\sum_{\ell=1}^d\dk_{\ell,v})$ is equal to either $0$ or $h$.
\end{enumerate}
\end{lemma}
\begin{proof}
First we claim that for any integer $k\geq 0$  we have
\begin{eqnarray*}
\llfloor\log_2(kd+1)\rrfloor \leq k\llfloor\log_2(d+1)\rrfloor,
\end{eqnarray*}
where the equality holds if and only if $k=0$ or $1$.
If $k=0$ or $1$ we see immediately that the equality holds.
If $k\geq 2$ then we have
\begin{eqnarray*}
\llfloor\log_2(kd+1)\rrfloor
&\leq& \llfloor \log_2 k + \log_2(d+1) \rrfloor\\
&\leq& (k-1)+\llfloor \log_2(d+1) \rrfloor\\
&<& (k-1)\llfloor \log_2(d+1) \rrfloor + \llfloor\log_2(d+1) \rrfloor\\
& = & k\llfloor \log_2(d+1)\rrfloor.
\end{eqnarray*}
This proves our assertion above.
Since $r=\sum_{v\geq 0}\sum_{\ell=1}^{d}\dk_{\ell,v}2^v$,
we have
\begin{eqnarray}\label{E:6}
s(r)&\leq& \sum_{v\geq 0}s(\sum_{\ell=1}^{d} \dk_{\ell,v})
\leq\sum_{v\geq 0}\llfloor\log_2(\sum_{\ell=1}^{d}\dk_{\ell,v}+1)\rrfloor
\\\nonumber
&\leq&\sum_{v\geq 0}\llfloor\log_2(\dk_{1,v}d+1)\rrfloor
\leq (\sum_{v\geq 0}\dk_{1,v}) \llfloor\log_2(d+1)\rrfloor
= s(\k)h.
\end{eqnarray}
Thus $s(\k)\geq \llceil\frac{s(r)}{h} \rrceil$.  Now suppose the
equality holds in (\ref{E:6}). Then $\dk_{1,v}=0$ or $1$. Because
$\dk_{\ell,v}\geq \dk_{\ell+1,v}$, Part 1) is immediate.  Since
$\sum_{v\geq 0}\dk_{1,v}=s(\k)=\frac{s(r)}{h},$ Part 2) follows.
Note that $s(\sum_{\ell=1}^{d}\dk_{\ell,v})\leq s(d)
\leq\llfloor\log_2(d+1)\rrfloor=h$. By Part 1), we know that
$s(r)$ is equal to the sum of
those $\frac{s(r)}{h}$ many nonzero $s(\sum_{\ell=1}^{d}\dk_{\ell,v})\leq
h$, thus each of the nonzero $s(\sum_{\ell=1}^{d}\dk_{\ell,v})$ has to
be equal to $h$. This proves 3).
\end{proof}

\begin{lemma}\label{L:binary}
Let notation be as above.
\begin{enumerate}
\item[a)]
For every $r\geq 1$ we have
$$\ord_2(C_r)\quad\geq\quad\llceil\frac{s(r)}{h}\rrceil.$$
\item[b)] For any integers $b, b'\geq 0$ we have
\begin{eqnarray*}
C_{2^{bh+b'}-2^{b'}}
&\equiv&C_{2^{bh}-1}^{2^{b'}} \bmod 2^{b+1}.
\end{eqnarray*}
\item[c)] If $g<2^h-2$ then
\begin{eqnarray*}
C_{2^{bh}-1}
&\equiv& 2^b(a_{2^{h}-1})^{\frac{2^{bh}-1}{2^{h}-1}}\bmod 2^{b+1}.
\end{eqnarray*}
\item[d)] If $g=2^h-2$ then
\begin{eqnarray*}
C_{2^{bh}-1}
&\equiv&2(a_{2^h-1})^{2^{(b-1)h}}C_{2^{(b-1)h}-1}\\
&&
+4(a_{3\cdot2^{h-1}-1})^{2^{(b-2)h}}C_{2^{(b-2)h}-1}\bmod 2^{b+1}.
\end{eqnarray*}
\end{enumerate}
\end{lemma}

\begin{proof}
a) By Lemma~\ref{miracle} and (\ref{E:bin}) we have
$\ord_2(C_r)\geq \min_{\k\in \K_r} {s(\k)}\geq \llceil \frac{s(r)}{h} \rrceil.$

b) For any natural number $t$ write $\gamma_t(\k)$ for the sum of the
entries in the $t$-th nonzero
column (from the left) of $\fbox{\k}$.

Let $r_b:=2^{bh}-1$. Suppose $\k\in\K_{2^{b'}r_b}$ with
$s(\k)=\frac{s(2^{b'}r_b)}{h}=b$. Then $\fbox{\k}$ consists of
only $0$ or $1$ entries and each
$s(\gamma_t(\k))=h$ by Lemma~\ref{miracle}. Since
$g\leq 2^h-2$, we have
$\gamma_t(\k)\leq 2(2^h-2)+1=2^{h+1}-3$. It follows that $(\gamma_t(\k))_2$
has at most $h+1$ digits, of which $h$ digits are equal to 1, and
at most one digit is equal to 0.
Moreover, this 0 can not be
at the rightest position. So the rightest $b'$ columns of $\fbox\k$
are zero. By removing those $b'$ zero columns of $\fbox\k$ we obtain
$\fbox{$\k'$}$ for $\k'\in\K_{r_b}$. It is clear that $s(\k')=b$.
Explicitly we have $k_\ell =2^{b'}k_\ell'$ for every $\ell =
1,\ldots,d$.  By~(\ref{E:bin}) we then have
\begin{eqnarray*}
C_{2^{b'}r_b}&\equiv& C_{r_b}^{2^{b'}}\bmod 2^{b+1}.
\end{eqnarray*}
This proves b).

c) Now let $\k\in\K_{r_b}$ with $s(\k)=b$. Arguing as above,
we observe that $\fbox{\k}$ consists of only $0$ or $1$ digits
and each $\gamma_t(\k)_2$ has
$h$ many 1's with at most one 0 in between. The 0 can only occur
if $g=2^h-2$. Indeed, if there is a 0 in
between, then either
$(\gamma_{t-1}(\k))_2$ or $(\gamma_{t+1}(\k))_2$
has to contribute a 1 at the corresponding position of $(r_b)_2$.
This situation can only happen if for $\mu=t$ or $\mu=t-1$ we have
$$
\begin{array}{l}
(\gamma_{\mu}(\k))_2
=(2^{h+1}-3)_2=
\overbrace{11\ldots11}^{h-1}01,\\
(\gamma_{\mu+1}(\k))_2
=(3\cdot2^{h-1}-1)_2=10\overbrace{11\ldots11}^{h-1}
\end{array}.
$$
However, $\gamma_{\mu}(\k)=2^{h+1}-3$ can only occur if $g=2^h-2$.

Suppose $g<2^h-2$. Then
$(\gamma_t(\k))_2
=2^h-1$ for $t=1,\ldots, b$,
and in between any two consecutive nonzero columns there are $h-1$ zero
columns.  It is then clear that $\k={}^t\!(k_1,\ldots,k_d)
\in\K_{r_b}$ is
defined by $k_1=k_2=\cdots =k_{2^{h}-1} =
\frac{r_b}{2^{h}-1}$ and
$k_{2^{h}}=\cdots=k_d=0$.  By $(\ref{E:bin})$
we have
$$C_{r_b}\equiv 2^{b}\left(a_{2^{h}-1}\right)^
{\frac{r_b}{2^{h}-1}} \bmod 2^{b+1}.$$
This proves c).

Suppose that $g=2^h-2$. We argued above that either
$\gamma_1(\k)=2^{h}-1$ or $\gamma_1(\k)=2^{h+1}-3$.  If
$\gamma_1(\k)=2^{h}-1$ then we remove the leftest non-zero
column from $\fbox{\k}$. This way, one obtains \fbox{$\k'$} for
$\k'\in\K_{r_{b-1}}$ with $s(\k')=b-1$. The relation between $\k$ and
$\k'$ is given by $k'_\ell=k_\ell-2^{(b-1)h}$ for
$\ell<2^{h}-1$ and $k'_\ell=k_\ell$ for $\ell\geq 2^{h}-1$. This
gives a $1-1$ correspondence between the two sets
$$\{\k\in\K_{r_b}\mid s(\k)=b{\rm\ and\ } \gamma_1(\k)=2^h-1\}$$ and
$$\{\k'\in\K_{r_{b-1}}\mid s(\k')=b-1\}.$$

If
$\gamma_1(\k)=2^{h+1}-3$ then we define $\k''\in\K_{r_{b-2}}$ with
$s(\k'')=b-2$ by removing the two leftest non-zero columns of
$\fbox\k$. In this case
$$k''_\ell=k_\ell-2^{(b-1)h-1}-2^{(b-2)h}$$ for $0\leq\ell<
3\cdot2^{h-1}-1$ and $k''_\ell=k_\ell-2^{(b-1)h-1}$ for $\ell\geq
3\cdot2^{h-1}-1$. This gives a $1-1$ correspondence between the two sets
$$\{\k\in\K_{r_b}\mid s(\k) = b{\rm\ and\ } \gamma_1(\k)=2^{h+1}-3\}$$ and
$$\{\k''\in\K_{r_{b-2}}\mid s(\k'') = b-2\}.$$ From $(\ref{E:bin})$ it now
easily follows that
\begin{eqnarray}
\label{telescope}
C_{r_b}&\equiv& 2\left(a_{2^h-1}\right)^{2^{(b-1)h}}
C_{r_{b-1}}+4\left(a_{3\cdot2^{h-1}-1}\right)^{2^{(b-2)h}}
C_{r_{b-2}}\ \bmod 2^{b+1}.
\end{eqnarray}
This proves d).
\end{proof}

\section{Proof of Theorem~\ref{MT:3}}
         \label{section3}

\begin{proof}[Proof of Theorem~\ref{MT:3}] From
Lemma~\ref{L:binary} it
follows that
$$\ord_2(C_{m2^{n+g-1}-j})\geq \llceil \frac{s(m2^{n+g-1}-j)}{h}\rrceil
\geq
\llceil\frac{n+g-h}{h}\rrceil\geq\llceil \frac{n}{h}\rrceil$$
for all $m,n \geq 1$.
By Key-Lemma~\ref{keylemma} i) we have
$\NP_1(X)\geq\frac{1}{h}$. This proves (I).
(See also~\cite{SZ:2}.)

{From} now on we assume $\NP_1(X)>\frac{1}{h}$.  For any
integer $n>1$ define $$\lambda_n:=\frac{n+g-2}{h(n-1)}.$$ Consider
$\lambda_n$ as a function in $n$, it is clear that $\lambda_n$ is
monotonically decreasing and converges to $\frac{1}{h}$ when $n$
approaches $\infty$.  Choose  $n_0$ such that
$\lambda_{n_0}<\NP_1(X)$ and such that $n_0+g-1$ is a multiple of
$h$ and such that $\frac{g-1}{h(n_0-1)}\leq 1$.
For all $1\leq n<n_0$ we have
$\lambda_{n_0}\leq\lambda_{n+1}=\frac{n+g-1}{hn}$; that is,
\begin{eqnarray*}
\llceil n\lambda_{n_0}\rrceil&\leq&\llceil\frac{n+g-1}{h}\rrceil.
\end{eqnarray*}
Therefore, for all $m\geq 1$ and $1\leq n<n_0$ one has
\begin{eqnarray*}
\ord_2(C_{m2^{n+g-1}-1})\geq \llceil \frac{n+g-1}{h}\rrceil
\geq \llceil n\lambda_{n_0}\rrceil.
\end{eqnarray*}
On the other hand, it follows that
\begin{eqnarray*}
\llceil n_0\lambda_{n_0}\rrceil
=\frac{n_0+g-1}{h}+\llceil\frac{g-1}{h(n_0-1)}\rrceil
=\frac{n_0+g-1}{h}+1.
\end{eqnarray*}
Hence, for all $m\geq 2$ one has
\begin{eqnarray*}
\ord_2(C_{m2^{n_0+g-1}-1})\geq\llceil\frac{n_0+g}{h}\rrceil
=\frac{n_0+g-1}{h}+1
=\llceil n_0\lambda_{n_0}\rrceil.
\end{eqnarray*}
Thus the hypotheses of Key-Lemma~\ref{keylemma} ii)
are satisfied. But $\NP_1(X)>\lambda_{n_0}$, so
we have
\begin{eqnarray}\label{E:bound1}
\ord_2(C_{2^{n_0+g-1}-1})\geq \llceil n_0\lambda_{n_0}\rrceil.
\end{eqnarray}

Now we apply our assumption that $g<2^h-2$ to
Lemma~\ref{L:binary} c) and get
\begin{eqnarray}\label{E:bound2}
C_{2^{n_0+g-1}-1} &\equiv& 2^{\llceil n_0\lambda_{n_0}\rrceil-1}
a_{2^{h}-1}^{(2^{n_0+g-1}-1)/(2^{h}-1)}\bmod
2^{\llceil n_0\lambda_{n_0}\rrceil}.
\end{eqnarray}
Hence $a_{2^{h+1}-1}\equiv 0\bmod 2$ by combining~(\ref{E:bound1})
and~(\ref{E:bound2}). Thus $c_{2^{h+1}-1}=0$.
This proves (II).

Assume $g=2^h-2$.
Define $$\lambda'_n := \frac{n+g-h-2}{h(n-2)}.$$ Consider
$\lambda'_n$ as a function in $n$, it is monotonically decreasing to
$\frac{1}{h}$ as $n$ approaches $\infty$.  Recall from above that
$n_0+g-1$ is a multiple of $h$ and $\lambda_{n_0}<\NP_1(X)$. Now we
may assume $n_0$ above is chosen large enough so that
$\lambda'_{n_0}<\NP_1(X)$ and $\frac{g-h}{h(n_0-1)}\leq 1$.
Since for all $n\leq n_0-2$ we have
$\lambda'_{n_0}\leq \lambda'_{n+2}$, we get
\begin{eqnarray*}
\llceil n\lambda'_{n_0}\rrceil&\leq&\llceil\frac{n+g-h}{h}\rrceil.
\end{eqnarray*}
It follows that for all $m\geq 1$ and $1\leq n\leq n_0-2$ one has
\begin{eqnarray*}
\ord_2(C_{m2^{n+g-1}-2^{h-1}})\geq\llceil\frac{s(2^{n+g-1}-2^{h-1})}{h}\rrceil=
\llceil\frac{n+g-h}{h}\rrceil
\geq\llceil n\lambda'_{n_0}\rrceil.
\end{eqnarray*}
On the other hand,
\begin{eqnarray*}
\llceil(n_0-1)\lambda'_{n_0}\rrceil
&=&\frac{n_0+g-h-1}{h}+\llceil\frac{g-h}{h(n_0-2)}\rrceil
\quad =\quad \frac{n_0+g-h-1}{h}+1.
\end{eqnarray*}
Therefore, for all $m\geq 2$ one has
\begin{eqnarray*}
\ord_2(C_{m2^{n_0+g-2}-2^{h-1}})\geq\llceil\frac{n_0+g-h}{h}\rrceil=
\frac{n_0+g-h-1}{h}+1=\llceil(n_0-1)\lambda'_{n_0}\rrceil.
\end{eqnarray*}
Hence the hypotheses of Key-Lemma~\ref{keylemma} ii) are satisfied.
But $\NP_1(X)>\lambda'_{n_0}$ implies that
$$\ord_2(C_{2^{n_0+g-2}-2^{h-1}})\quad\geq\quad \llceil
(n_0-1)\lambda'_{n_0}\rrceil.$$
By Lemma~\ref{L:binary}~b), this implies
\begin{eqnarray}\label{Ezero}
C_{2^{n_0+g-h-1}-1}&\equiv& 0 \bmod 2^{\llceil(n_0-1)\lambda'_{n_0}\rrceil}.
\end{eqnarray}
Recall from~(\ref{E:bound1}) we have
$$\ord_2(C_{2^{n_0+g-1}-1})\geq \llceil n_0\lambda_{n_0}\rrceil
=\frac{n_0+g-1}{h}+1.$$
Put these together we get for $b=\frac{n_0+g-1}{h}$ that
\begin{eqnarray}\label{E:telescope}
C_{2^{bh-1}}\equiv 2C_{2^{(b-1)h}-1}\equiv 0\bmod
2^{b+1}.
\end{eqnarray}

Now recall
Lemma~\ref{L:binary}~d), for $3\leq b\leq \frac{n_0+g-1}{h}$,
we have
\begin{eqnarray}\label{Enonzero}
C_{2^{bh}-1} &\equiv& 2(a_{2^h-1})^{2^{(b-1)h}}
C_{2^{(b-1)h}-1}\\\nonumber
&&+4(a_{3\cdot2^{h-1}-1})^{2^{(b-2)h}}C_{2^{(b-2)h}-1}
\bmod 2^{b+1}.
\end{eqnarray}
Assume
$a_{3\cdot2^h-1}\not\equiv 0\bmod 2$ or $a_{2^h-1}\not\equiv 0
\bmod 2$. By applying~(\ref{Enonzero}) recursively for
$b=\frac{n_0+g-1}{h},\ldots, 4,3$,
and also by using (\ref{E:telescope}) one gets for all $b\geq1$ that
$$C_{2^{bh}-1}\equiv 0 \bmod 2^{b+1}.$$
When $b=1$
we have $$2a_{2^{h}-1}\equiv
C_{2^{h}-1}\equiv0\bmod 2^2$$ it follows that $c_{2^{h}-1}=0$.
Apply this to the congruence of $b=2$ one get
$$4a_{3\cdot2^{h-1}-1}\equiv C_{2^{2h}-1}\equiv0\bmod 2^3$$
hence $c_{3\cdot2^{h-1}-1}=0$.  This proves (III).
\end{proof}

\section{Moduli of hyperelliptic supersingular curves}
         \label{section4}

The goal of this section is
to prove Theorems~\ref{MT:4} and~\ref{MT:2}, as an application of
Theorem~\ref{MT:3}.

\begin{proposition}\label{P:1}
Let $X$ be a hyperelliptic curve over $\bar\F_2$ of genus $g$ and
$2$-rank zero. Let $m< 2g$ be an odd integer, such that the binomial
coefficient $\binom{2g+1}{2^km}$ is odd for some $k\geq 0$.  Then
there exists an equation for $X$ of the form
\begin{eqnarray}\label{E:oddterms}
y^2-y=\sum_{i=0}^{g}c_{2i+1}x^{2i+1}
\end{eqnarray}
with $c_{2g+1}=1$ and $c_m=0$. Moreover, there are only finitely
many equations for $X$ of this form.
\end{proposition}
\begin{proof}
Every $2$-rank zero hyperelliptic curve of genus $g$ has an equation
of the form $y^2-y=\tf(x)$ where $\tf(x)$ is a polynomial of degree $d$.
This follows from the Deuring-Shafarevich formula
(see~\cite{Nakajima:85}) as an easy exercise. By
applying an obvious isomorphism of $X$ we can assume $\tf$ is monic and
odd. From now on we write $\tf(x)=\sum_{i=0}^gc_{2i+1}x^{2i+1}$ with
$c_{2g+1}=1$.
Any isomorphism between curves with an equation of the form
$(\ref{E:original})$ is of the form
$\sigma: (x,y)\mapsto (\zeta x+t_0,y+h(x))$ for some
polynomial $h(x)$ of degree $\leq g$, some
$t_0\in\bar\F_2$ and some $d$-th
root of unity $\zeta$. (This polynomial $h(x)$ has nothing to
do with the integer $h$ defined in this paper.)
We
want to show that there exists
a triple $(t_0,\zeta,h(x))$ such that $\sigma$ sends
$X$ to
$y^2-y=\tf(\zeta x+t_0)+h(x)^2+h(x)$,
where the polynomial on the right-hand-side is odd in $x$ and
has its $x^m$-coefficient
vanishing. Also, we want to show that
there are only finitely many such $(t_0,\zeta,h(x))$.
Let $t$ be a variable. Write out
$\tf(\zeta x+t)=\sum_{j=0}^{2g+1}f_j(t)x^j$ for polynomials $f_j(t)\in
\bar\F_2[t]$. Let $h(x)=\sum_{i=0}^g h_i(t)x^i$ for some
$h_i(t)\in\bar{\F_2(t)}$ such that
$\tf^\sigma:= \tf(\zeta x+t)+h(x)^2+h(x)$ is odd in $x$.
Let $h_i(x)=0$ for $i>g$.
It suffices to show that the
$x^m$-coefficient $f_m+h_m$ of $\tf^\sigma$ is a non-constant function in
$t$ so that it vanishes after specializing $t$ to one of finitely many
$t_0\in\bar\F_2$.

Since $\tf^\sigma(x)$ is odd, for $i>g/2$ we have
$h_i=\sqrt{f_{2i}}$, and for $0<i\leq g/2$ we have recursively
$h_i=\sqrt{f_{2i}+h_{2i}}$. For $j$ even, $f_j(t)$ is odd with
a term $\binom{d}{j}t^{d-j}$.
So by the hypothesis on $m$, there exists a $k$ such that
$f_{2^km}\not=0$. Let $s$ be the biggest $k$ with
$f_{2^km}\not=0$.
Then we have $h_{2^km}=0$ for all $k\geq s$.
We consider the cases $s>0$ and $s=0$
separately.

If $s>0$ then $h_{2^{s-1}m}=\sqrt{f_{2^{s}m}}$.
Recall that $f_{2^{s}m}$ is odd, hence $h_{2^{s-1}m}\in
{\overline{\F_2(t)}}\backslash\bar\F_2(t)$.
For any $A\in{\overline{\F_2(t)}}\backslash\bar\F_2(t)$ one has that
$\sqrt{A}\in{\overline{\F_2(t)}}\backslash\bar\F_2(t)$. By applying this
$s-1$ times one observes that $h_m\in
{\overline{\F_2(t)}}\backslash\bar\F_2(t)$.

If $s=0$ then $\binom{d}{2^km}$ is even for all $k>s$.
Hence $\binom{d}{m}$ is odd, $f_m=\binom{d}{m}t^{d-m}$ and
$h_m=0$.

In both cases $f_m+h_m$ is not constant. This finishes the proof.
\end{proof}

\begin{proof}[Proof of Theorem~\ref{MT:4}]
By Theorem~\ref{MT:3} and Proposition~\ref{P:1}, any supersingular
hyperelliptic curve of genus $g$
can be given by an equation of
the form $(\ref{E:oddterms})$ with $c_1=c_{2^{h}-1}=0$. If
$g=2^{h}-2$ then we also have $c_{3\cdot2^{h-1}-1}=0$.
For $g<2^{h}-2$ (respectively, $g=2^{h}-2$)
there are only $g-2$ (respectively, $g-3$) non-zero coefficients
left.
\end{proof}

Below we exhibit $2$-parameter families of hyperelliptic
supersingular curves of genus $4$. So its number of moduli
is equal to $2$ by Proposition~\ref{P:1} and hence Theorem~\ref{MT:2}
follows immediately.

\begin{proof}[Proof of Theorem~\ref{MT:2}]
The given curve is of the form in~(\ref{E:geer}) hence supersingular.  By
Proposition~\ref{P:1} every genus-$4$ hyperelliptic curve of $2$-rank
zero has an equation $y^2-y = x^9+c_7x^7+c_5x^5+c_3x^3$ for some
$c_i\in\bar\F_2$, which is not supersingular if $c_7\neq 0$ by
Theorem~\ref{MT:3}.
\end{proof}

\begin{remark}
Normal forms of hyperelliptic supersingular curves
of low genera are studied in more details in an upcoming
paper~\cite{SZ:3}.
\end{remark}

\end{document}